\documentclass[a4paper,english]{scrartcl}
\usepackage[utf8x]{inputenc}
\usepackage[T1]{fontenc}
\usepackage{babel}
\usepackage{amsmath,amssymb}
\usepackage[amsmath,thref,thmmarks]{ntheorem}
\usepackage[pdftex]{graphicx}
\usepackage{color}
\usepackage{url}
\usepackage[all]{xy}

\theoremstyle{plain}
\theoremheaderfont{\normalfont\bfseries}
\theorembodyfont{\itshape}
\theoremseparator{:}
\theoremnumbering{arabic}
\theoremsymbol{}

\newtheorem{thm}{Theorem}[section]

\theoremstyle{plain}
\theoremheaderfont{\normalfont\bfseries}
\theorembodyfont{\itshape}
\theoremseparator{:}
\theoremnumbering{arabic}
\theoremsymbol{}

\newtheorem{proposition}[thm]{Proposition}

\theoremstyle{plain}
\theoremheaderfont{\normalfont\bfseries}
\theorembodyfont{\itshape}
\theoremseparator{:}
\theoremnumbering{arabic}
\theoremsymbol{}

\theoremstyle{plain}
\theoremheaderfont{\normalfont\bfseries}
\theorembodyfont{\itshape}
\theoremseparator{:}
\theoremnumbering{arabic}
\theoremsymbol{}

\newtheorem{lemma}[thm]{Lemma}

\theoremstyle{plain}
\theoremheaderfont{\normalfont\bfseries}
\theorembodyfont{\itshape}
\theoremseparator{:}
\theoremnumbering{arabic}
\theoremsymbol{}

\theoremstyle{plain}
\theoremheaderfont{\normalfont\bfseries}
\theorembodyfont{\itshape}
\theoremseparator{:}
\theoremnumbering{Alph}
\theoremsymbol{}

\theoremstyle{plain}
\theoremheaderfont{\normalfont\bfseries}
\theorembodyfont{\itshape}
\theoremseparator{:}
\theoremnumbering{Alph}
\theoremsymbol{}

\theoremstyle{plain}
\theoremheaderfont{\bfseries}
\theorembodyfont{\normalfont}
\theoremseparator{:}
\theoremnumbering{arabic}
\theoremsymbol{}

\newtheorem{definition}[thm]{Definition}

\theoremstyle{plain}
\theoremheaderfont{\bfseries}
\theorembodyfont{\normalfont}
\theoremseparator{:}
\theoremnumbering{arabic}
\theoremsymbol{}

\theoremstyle{plain}
\theoremheaderfont{\bfseries}
\theorembodyfont{\normalfont}
\theoremseparator{:}
\theoremnumbering{arabic}
\theoremsymbol{}

\theoremstyle{nonumberplain}
\theoremheaderfont{\bfseries}
\theorembodyfont{\normalfont}
\theoremseparator{:}
\theoremnumbering{arabic}
\theoremsymbol{}

\newtheorem{bemon}{Remark}

\theoremstyle{nonumberplain}
\theoremheaderfont{\bfseries}
\theorembodyfont{\normalfont}
\theoremseparator{:}
\theoremnumbering{arabic}
\theoremsymbol{}

\theoremstyle{plain}
\theoremheaderfont{\normalfont\scshape}
\theorembodyfont{\itshape}
\theoremseparator{:}
\theoremnumbering{arabic}
\theoremsymbol{}

\newcounter{QuesCounter}
\theoremstyle{plain}
\theoremheaderfont{\normalfont\scshape}
\theorembodyfont{\itshape}
\theoremseparator{:}
\theoremnumbering{arabic}
\theoremsymbol{}

\newtheorem{frage}[QuesCounter]{Question}

\theoremstyle{nonumberplain}
\theoremheaderfont{\normalfont\scshape}
\theorembodyfont{\itshape}
\theoremseparator{:}
\theoremnumbering{arabic}
\theoremsymbol{}

\theoremstyle{nonumberplain}
\theoremseparator{:}
\theoremheaderfont{\scshape}
\theorembodyfont{\upshape}
\theoremsymbol{{\Large\ensuremath{_\dashv}}}

\newtheorem{beweis}{Proof}

\theoremstyle{nonumberplain}
\theoremseparator{:}
\theoremheaderfont{\scshape}
\theorembodyfont{\upshape}
\theoremsymbol{\ensuremath{\square}}

\newcommand{\<}{\langle}
\renewcommand{\>}{\rangle}
\renewcommand{\phi}{\varphi}

\newcommand{\gdw}{\leftrightarrow}
\newcommand{\op}[1]{\operatorname{#1}}
\newcommand{\Pot}{\op{\mathcal{P}}}

\newcommand{\union}{\bigcup\limits}
\newcommand{\Union}{\bigcup}

\newcommand{\restr}{\upharpoonright}

\newcommand{\ult}{\op{Ult}}
\newcommand{\on}{\op{On}}

\newcommand{\crit}{\op{crit}}
\newcommand{\ran}{\op{ran}}

\newcommand{\card}{\op{Card}}
\newcommand{\cof}{\op{cof}}

\newcommand{\Sk}{\op{Sk}}
\newcommand{\hull}{\op{Hull}}

\newcommand{\ptwimg}[2]{{#1}"\left[{#2}\right]}
\newcommand{\kleiner}{\mathord{<}}

\newcommand{\mitord}{\op{o}}

\newcommand{\finsubsets}[1]{{\left[#1\right]}^{\kleiner\omega}}

\newcommand{\eextend}{\trianglelefteq}
\newcommand{\propeextend}{\triangleleft}

\begin{document}

\title{\bfseries Lower consistency bounds for mutual stationarity with divergent cofinalities and limited covering}
\author{Dominik Adolf}
\maketitle

\begin{abstract}
	We improve previous work on the consistency strength of mutually stationary sequences of sets concentrating on points with divergent cofinality building on previous work by Adolf, Cox and Welch. Specifically, we have greatly reduced our reliance on covering properties in the proof. This will allow us to handle sequences in which sets concentrating on points of countable cofinality appear infinitely often. Furthermore we will show that if $\kappa$ is a J\'{o}nsson cardinal with $\kappa < \aleph_\kappa$ then $0^\P$, the sharp for a model with a strong cardinal, exists. 
\end{abstract}

\section{Introduction}

\begin{definition}[Foreman-Magidor]
	Let $\<\kappa_i: i < \delta\>$ be a strictly increasing sequence of regular uncountable cardinals. A sequence $\<S_i:i < \delta\>$ with $S_i \subseteq \kappa_i$ for all $i < \delta$ is \textit{mutually stationary} if and only if the set of $A \subseteq \lambda := \sup\limits_{i < \delta} \kappa_i$ with $\sup(A \cap \kappa_i) \in S_i$ whenever $\kappa_i \in A$ is stationary, i.e. for all $F:\finsubsets{\lambda} \rightarrow \lambda$ there is some such $A$ with $\ptwimg{F}{\finsubsets{A}} \subseteq A$.
\end{definition}

The property of mutual stationarity was first introduced by Foreman and Magidor in \cite{mutstat}. In that paper they proved that any sequence of stationary sets concentrating on points of countable cofinality is mutually stationary. They also proved that the same does not hold in $L$ with sequences concentrating on points of cofinality $\omega_1$. ( By work of Koepke and Welch this property does have large cardinal strength \cite{strengthmutstat}. Ben-Neria has recently shown its consistency in \cite{singularstatI}. It seems likely that it is quite strong.)

Here though we will be solely interested in sequences that are not limited to a single cofinality. In fact, we do consider sequences with certain recurring patterns. Improving on \cite{altmutstat} we show:  

\begin{thm}\label{first}
 Assume that $0^\P$ does not exist. Let $2 \leq k,l < \omega$. Let $\<S_n : k \leq n < \omega\>$ be a sequence such that
 \begin{itemize}
  \item $S_n \subset \aleph_n$ is a stationary set that concentrates on limit ordinals with a fixed cofinality $\mu_n$;
  \item the sequence $\<S_n: k \leq n < \omega\>$ is mutually stationary;
  \item there exists a sequence $\<n_i: i < \omega\>$ such that 
        \begin{itemize}
         \item $n_{i + 1} \geq n_i + l$,
         \item $\mu_{n_i} = \mu_{n_i + j}$ for $j < l$,
         \item more than one value appears infinitely often in the sequence $\<\mu_{n_i}: i < \omega\>$.
        \end{itemize}
 \end{itemize}
 Then there exist infinitely many $n < \omega$ such that there exists $\kappa_n < \aleph_n$ such that $(\kappa^+_n)^K < \aleph_n$, $\mitord^K(\kappa_n) \geq \max(\aleph_n,(\kappa^{+(l+1)}_n)^K)$.
\end{thm}

Important here is that we do allow the $\mu_n$ to be $\omega$ infinitely often. Unlike with Theorem 6 of \cite{altmutstat} where cofinalities track closely with the Mitchell order of measures on $K$ the actual values of the $\mu_n$ seemed of no consequence in \cite[Theorem 7]{altmutstat}, only the pattern in which they appear. Of course, we did assume there that all the $\mu_n$ were uncountable, but the new result indicates that this requirement was not necessary.

By work of Shelah \cite{liushelah} and Sharon \cite{sharonmsth} the existence of a sequence as above is consistent relative to the existence of a cardinal $\kappa$ with $\mitord(\kappa)\geq \kappa^{+\omega}$. Ben-Neria has suggested that a careful reading of the argument yields that in this specific case this bound can be reduced down to the existence of a cardinal $\kappa$ with $\mitord(\kappa) \geq \kappa^{+(l+1)}+1$. The sequence of $\<\kappa_n: n < \omega\>$ given by the theorem are not fully unlike a Prikry sequence which suggests that this might indeed be an exact bound.

%In the paper we prove that in $K$ there exist unboundedly in $\aleph_\omega$ many $\kappa$ with $\mitord(\kappa) \geq \kappa^{+(l+1)}$. A quick observation will yield significantly more strength.

%By $(\mitord(\kappa) \geq \kappa^{+(l+1)})^\ddag$ we refer to the least sound active $(\omega_1 + 1)$-iterable premouse $(M;\in,\vec{E},F)$ such that $M \vert\vert\crit(F)$ believes there are proper class many $\kappa$ with $\mitord(\kappa)\geq \kappa^{+(l+1)}$.

Improving \cite[Theorem 7]{altmutstat} also allows us to improve \cite[Theorem 8]{altmutstat}:

\begin{thm}\label{second}
 Let $k < l$ be natural numbers. Assume that for all $f:\omega \rightarrow \{k,l\}$ the sequence $\<S^n_{f(n)}: k < n < \omega\>$ is mutually stationary. Then $0^\P$ exists.
\end{thm}

Again the difference is that we do allow countable cofinalities. This is significant, because for $k = 0$ this property is known to be consistent by work of Shelah \cite[Section 6]{canstrtwo}. (For $l = 1$ Jensen has a consistency proof with a lesser requirement.(\cite{chforcax}))

Our last theorem improves on \cite[Theorem 9]{altmutstat} in which we proved that the existence of $0^\P$ follows from the existence of a mutually stationary sequence $\<S_n:n < \omega\>$ concentrating on cofinalities $\<\mu_n:n < \omega\>$ such that no cofinality appears more than finitely often among the $\<\mu_n: n < \omega\>$. The existence of such a sequence implies that $\aleph_\omega$ is J\'{o}nsson, and our improved theorem will not assume more than that.

This is a departure from the arguments we used previously: while a typical J\'{o}nsson structure $X \subset \aleph_\omega$ will have $\<k_n: n < \omega\>$ with $\cof(X \cap \aleph_{k_n}) = \aleph_n$ for all $n$, this is not quite enough control about all the possible cofinalities that can appear during the iteration. We get around this by using the fact that many levels of the J\'{o}nsson structure are instances of Chang's conjecture together with using a "pseudo-drop" used by Mitchell in \cite{mitchelljonsson} to fix cofinalites on a club set. A contradiction can then be reached using our usual methods.

\begin{thm}\label{third}
 Assume that $\kappa$ is J\'{o}nsson with $\kappa < \aleph_\kappa$. Then $0^\P$ exists.
\end{thm}

The rest of this paper will be organized in the following fashion: the next section will review the necessary facts about Inner Model Theory and stationary sets used in the proof and introduce some useful definitions; section 3 will contain the proofs of \thref{first} and \thref{second}; section 4 will contain the proof of \thref{third}; we will finish with open questions and acknowledgements.

\section{Preliminaries}

In the following we will often confuse some $S \subset \Pot(\kappa)$ with $S^+ := \{ A \subseteq H_\kappa \vert A \cap \kappa \in S\}$. It is a standard fact that this makes no difference as far as stationarity is concerned. 

We will use the same approach to fine structure as in the previous paper \cite{altmutstat}, so we use the general outline of \cite{Zeman}. $K$ where it appears will always refer to the core model below $0^\P$.

Let $X \prec H_\kappa$ for some uncountable cardinal $\kappa$. We write $H_X$ for the transitive collapse of $X$, and $\sigma_X: H_X \rightarrow X$ for the unique isomorphism. As usual we will want to co-iterate $K$ and $K_X := \ptwimg{\sigma^{-1}_X}{K}$. (What structure specifically $X$ is a substructure of will vary between sections, but will be clear from context.)

When co-iterating $K$ and $K_X$ let $\theta_X$ be the length of the iteration on the $K$-side, let $M^X_i$ ($i \leq \theta_X$) be the $i$-th model appearing in that iteration, $\kappa^X_i$ the critical point of the $i$-th extender used during the iteration, and $\nu^X_i$ the corresponding iteration index. We let $\pi^X_{i,j}:M^X_i \rightarrow M^X_j$ ($i \leq j \leq \theta_X$) be the iteration embedding(, if there is a truncation in the interval $\left(i,j\right]$ then $\pi^X_{i,j}$ will be partial). Finally, we will let $d^X_i$ be the eventual degree of elementarity of the embeddings into $M^X_i$. (Note that if $i \leq j \leq \theta_X$, there is no truncation in the interval $\left(i,j\right]$ and $d^X_\cdot$ is constant there, then $\pi^X_{i,j}$ is continuous at $\rho_{d^X_i}(M^X_i)$.)

In Section 5 we will have to consider special iterations on $K$ in which we might be forced to use a special form of truncation. We shall still use the same notation for these special iterations. This should not lead to confusion.

Unlike in \cite{altmutstat}, in the arguments to come we will not assume that the $K_X$-side remains trivial in these co-iterations both special and regular. We will write $\zeta_X$ for the length of the iteration on the $K_X$-side. We will write $N^X_i$ ($i \leq \zeta_X$) for  the $i$-th model appearing in that iteration. For $i \leq j \leq \zeta_X$ we will write $\tau^X_{i,j}: N^X_i \rightarrow N^X_j$ for the iteration embedding. (Note: There will be no truncations on the $K_X$-side.)

We will make heavy use of the following lemma (, this is Lemma 10 in \cite{altmutstat}):

\begin{lemma}\label{corelemma}
 Let $M$ be a $J$-structure and $n < \omega$. Let $\lambda,\kappa$ be cardinals in $M$,  $\kappa$ regular in $M$, $\rho_{n+1}(M) \leq \lambda < \kappa \leq \rho_{n}(M)$, and $M$ is $n+1$-sound above $\lambda$, i.e. $\hull^M_{n + 1}(\lambda \cup \{p^M_{n+1}\}) = M$. Then $\cof(\kappa) = \cof(\rho_n(M))$.
\end{lemma}

Note here $\hull^M_{n+1}$ is a fine-structural Skolem hull, i.e. relative to the canonical $\Sigma_1$ Skolem functions over the $n$-th reduct of $M$. By the nature of stationary sets we will have to consider Skolemized structures $\mathfrak{A}$ on $H_\kappa$. In that case we will write $\Sk^\mathfrak{A}(A)$ for the appropriate Skolem hull of $A \subseteq H_\kappa$.

\section{Mutually stationary sequences with alternating blocks}

We start out by introducing and subsequently analysing some terms and definitions that will serve to better communicate the argument.

\begin{definition}
	Let $S \subset \Pot(\aleph_\omega)$ and $2 < l < \omega$.
	\begin{itemize}
		\item[$(a)$] $S$ has fixed cofinalities (with values $\<\mu_n: n < \omega\>$) if and only if $\cof(X \cap \aleph_n) = \mu_n$ for all $X \in S$.
		\item[$(b)$] $S$ has alternating blocks of size $l$ if and only if $S$ has fixed cofinalities with values $\<\mu_n: n < \omega\>$ and there exists $\<n_k: k < \omega\>$ such that $\mu_{n_k} = \mu_{n_k + i}$ for all $i < l$ and more than one value appears infinitely often in the sequence $\<\mu_{n_k}:k < \omega\>$. 
	\end{itemize}
\end{definition}

The hypothesis of \thref{first} (relative to $2 < l < \omega$) can thus be restated to: "there exists a stationary $S \subset \Pot(\aleph_\omega)$ that has alternating blocks of size $l$". To increase distinctiveness between this section and the next we will make sure that our stationary set does not also witness the J\'{o}nsson-ness of $\aleph_\omega$. Furthermore we will see that having alternating blocks of any size is stable under changes done to an initial segment of our structures.

\begin{proposition}\label{prop1}
	Let $S \subset \Pot(\aleph_\omega)$ and $2 < l < \omega$.
	\begin{itemize}
		\item[$(a)$] If $S$ has alternating blocks of size $l$ and is stationary, then there exists some $q < \omega$ and some stationary $S^* \subset \Pot(\aleph_\omega)$ with alternating blocks of size $l$ such that $\card(X) <\aleph_q$ for all $X \in S^*$.
		\item[$(b)$] Let $k < \omega$. If $S$ has alternating blocks of size $l$, is stationary, and there exists $q < \omega$ with $\card(X) < \aleph_q$ for all $X$ in $S$, then there exists some stationary $S^* \subset \Pot(\aleph_\omega)$ with alternating blocks of size $l$ such that $\card(X) < \aleph_q$ and $\cof(X \cap \aleph_n) \geq \aleph_1$ for all $X \in S^*$ and $0 < n < k$.
	\end{itemize}
\end{proposition}

\begin{beweis}
	Let $S \subset \Pot(\aleph_\omega)$ stationary and $2 < k,l < \omega$. Assume that $S$ has blocks of size $l$. We will first refine $S$ into some $S^*$ with uniformely bounded elements as required by $(a)$. We will then further refine to some stationary $S^{**}$ that fulfills the requirements of $(b)$.
	
	Let $\<\mu_n: n < \omega\>$ be the values associated to $S$ and $\<n_k : k  < \omega\>$ the starting points of each block. Let $\aleph_c,\aleph_d$ be two distinct values that appear infinitely often in the sequence $\<\mu_{n_k}:k < \omega\>$. We will see that $q := \max\{c,d\} + 1$ works.
	
	Let $\mathfrak{A}$ be some Skolemized structure on $\aleph_\omega$. By assumption we have some $X_\mathfrak{A} \in S$ with $X \prec \mathfrak{A}$. Let $C := \{k < \omega \vert \mu_{n_k} = \aleph_c\}$ and $D:= \{k < \omega \vert \mu_{n_k} = \aleph_d\}$. Fix $f_{k,i}: \mu_{n_k + i} \rightarrow X \cap \aleph_{n_k + i}$ cofinal. Let then $X^*_\mathfrak{A} : = \Sk^\mathfrak{A}(\Union\{\ran(f_{k,i})\vert k \in C \cup D, i < l\})$. It is then easy to check that $S^*$, the set of the $X^*_\mathfrak{A}$'s is as desired.
	
	From the $X^*_\mathfrak{A}$'s we can then construct a sequence $\<X^\alpha_\mathfrak{A}: \alpha < \omega_1\>$ with the following properties:
	\begin{itemize}
		\item $X^\alpha_\mathfrak{A} \prec \mathfrak{A}$ for all $\alpha < \omega_1$;
		\item $X^*_\mathfrak{A} \subseteq X^\alpha_\mathfrak{A} \subseteq X^\beta_\mathfrak{A}$ for all $\alpha \leq \beta < \omega_1$;
		\item $\sup(X^\alpha_\mathfrak{A} \cap \aleph_n) = \sup(X^*_\mathfrak{A} \cap \aleph_n)$ for all $\alpha < \omega_1$ and $n \geq k$;
		\item $\sup(X^\alpha_\mathfrak{A} \cap \aleph_n) < \sup(X ^\beta_\mathfrak{A} \cap \aleph_n)$ for all $\alpha < \beta < \omega_1$ and $n < k$ (assuming $\aleph_n \not\subset X^\alpha_\mathfrak{A}$).
	\end{itemize}
	A sequence like this is easily constructed: set $X^0_\mathfrak{A} := X^*_\mathfrak{A}$, $X^{\alpha + 1} := \Sk^\mathfrak{A}(X^\alpha_\mathfrak{A} \cup \{\sup(X^\alpha_\mathfrak{A} \cap \aleph_n):n < k\})$, and take unions at limit steps. That the third condition is fulfilled is due to a well-known lemma of Baumgartner's \cite{thatonelemma}.
	
	We then set $X^{**}_\mathfrak{A} := \union_{\alpha < \omega_1} X^\alpha_\mathfrak{A}$. $S^{**}$ the set consisting of the $X^{**}_\mathfrak{A}$'s is then as desired.
\end{beweis}

\begin{bemon}
	Note that we could have put $\aleph_{q - 1}$ into our hull without issue thanks to Baumgartner's lemma.
\end{bemon}

Let $2 < l < \omega$. From now on we shall assume that $0^\P$ does not exist. We thus have a core model $K$. Assume there exists some stationary set on $\Pot(\aleph_\omega)$ with alternating blocks of size $l$. As usual we will require that $K$ exhibits certain nice behaviours in co-iterations with hulls of itself. The next lemma will show that we will have many good hulls that have alternating blocks of size $l$.

\begin{lemma}
 There exists some $S \subset \Pot(H_{\aleph_\omega})$ stationary with alternating blocks of size $l$ such that $K$ truncates in co-iteration of itself with $K_X$ for all $X \in S$ (assuming the iteration lasts longer than $2$ steps).
\end{lemma}

\begin{beweis}
 Let $S$ be any stationary set with alternating blocks of size $l$. By \thref{prop1} we can assume that $\card(X) < \aleph_{q + 1}$, $\aleph_q \subset X$, and $\cof(X \cap \aleph_n) \geq \aleph_1$ for all $X \in S$ and $0 < n \leq q + 2$. We will follow the proof of (\cite{Cox}, Lemma 39) and adapt where necessary.
 
 Consider the co-iteration between $K_X$ and $K$. We do have that $K_X$ and $K$ agree up to $\alpha^X_0 := (\crit(\sigma_X)^+)^{K_X}$ which by assumption has uncountable cofinality. It follows that $\alpha^X_0$ cannot be a cardinal in $K$ as otherwise by \cite[Corollary 19]{Cox}, $\ult(K; \sigma_X \restr (K \vert\vert \alpha^X_0))$ is iterable.
 
 In conclusion the first extender applied to $K$ has index larger than $\alpha^X_0$. There are two cases: in the first case we have a stationary set $S' \subset S$ such that $\kappa^X_0 \geq \crit(\sigma_X)$ for all $X \in S'$. In that case the first extender applied to $K$ is not total over it and therefore $S'$ is as wanted.
 
 So, assume that $\kappa^X_0 < \crit(\sigma_X)$ for (almost) all $X \in S$. By pressing down we can assume that $\kappa^X_0$ has a fixed value $\kappa_0$. Then $\nu^X_0$ is the Mitchell-order of $\kappa_0$ in $K_X$. Then $\sigma_X(\nu^X_0)$ is $\mitord^K(\kappa_0) > \aleph_{q+1}$.
 
 Fix $r$ such that $(\mitord^K(\kappa_0)^+)^K < \aleph_r$. We can then find a stationary set $S' \subset \Pot(H_{\aleph_\omega})$ such that 
 \begin{itemize}
  \item $\card(X) < \aleph_{q+1}$, $\aleph_q \subset X$ for all $X \in S'$;
  \item $\crit(\sigma_X) > \kappa_0$ for all $X \in S'$;
  \item $\cof(X \cap \aleph_n) \geq \aleph_1$ for all $0 < n \leq r$;
 \end{itemize}
 
 and $S'$ has alternating blocks of size $l$. Let us now fix $X \in S'$. As before $\nu^X_0 \geq \alpha^X_0 \geq ((\kappa_0)^+)^K$. On the other hand $\nu^X_0 < \on \cap K_X < \aleph_{q + 1} < \mitord^K(\kappa_0)$. So the first extender applied to $K$ has critical point $\kappa_0$. 
 
 Thus $\nu^X_0$ is in fact $\mitord^{K_X}(\kappa_0)$. The local successor $\alpha^X_1 := ((\nu^X_0)^+)^{K_X}$ of the latter is less than $\sigma^{-1}_X(\aleph_r)$. (This is just the reflection of a first order fact about $K$ which we established previously.)
 
 By choice of $X$ we then have that $\alpha^X_1$ has uncountable cofinality (using the weak covering of $K$ reflected down to $K_X$). We conclude that $\alpha^X_1$ is not a cardinal in $M^X_1$ as otherwise $\ult(M^X_1; \sigma_X \restr (M^X_1 \vert\vert \alpha^X_1))$ is iterable. This means that $\alpha^X_1 < ((\nu^X_0)^+)^{M^X_1}$ and therefore the extender applied to $M^X_1$ is not total over it. Therefore $S'$ is as desired.
\end{beweis}

From now on fix some $S$ with alternating blocks of size $l$ as given by the lemma. Let $\<\mu_n: n < \omega\>$ be the cofinalities associated to $S$ and $\<n_k: k < \omega\>$ be a list of starting points for each block. For some $X \in S$ we let $\beta^X_n := \sigma^{-1}_X(\aleph_n)$. Obviously, $\cof(\beta^X_n) = \mu_n$ for all $n < \omega$.
 
It is not clear that we can show that infinitely many of the intervals $\left[\beta^X_{n_k},\beta^X_{n_k + (l-1)}\right]$ are overlapped by extenders appearing in the iteration. We will therefore make a slight substitution. Let $(\gamma^X_{n_i + j})^* := ((\beta^X_{n_i})^{+j})^{K_X}$ for $j < l$. Note that we still have $\cof((\gamma^X_{n_k + j})^*) = \mu_{n_k + j} (= \mu_{n_k})$. This should be clear in case $(\gamma^X_{n_k + j})^* = \beta^X_{n_k + j}$, but even if not it holds as a consequence of weak covering reflected down to $K_X$ (cf. \cite[Obs. 25]{altmutstat}). We then let $\gamma^X_{n_i + j} := \tau^X_{0,\zeta_x}((\gamma^X_{n_i + j})^*)$. Note:
 
\begin{enumerate}
  \item $\gamma^X_{n_i + j} = ((\gamma^X_{n_i})^{+j})^{N^X_{\zeta_X}}$ for $j < l$;
  \item $\cof(\gamma^X_{n_i + j}) = \mu_{n_i}$ for $0 < j < l$ (but not necessarily for $j = 0$ if $\beta^X_{n_k}$ is a measurable cardinal in $K_X$) because of the continuity of iteration embeddings at points of non-measurable cofinality.
\end{enumerate}

Now we can show:

\begin{lemma}
 $\theta_X$ is a limit ordinal for all $X \in S$. In fact, the generators of the iteration (on the $K$-side) are unbounded in $\on \cap N^X_{\zeta_X}$.
\end{lemma}

\begin{beweis}
 Assume not. As the $K$-side truncates we do have that $N^X_{\zeta_X} \propeextend M^X_{\theta_X}$. In $N^X_{\zeta^X}$, $\<\gamma^X_{n_i + 1} : i < \omega\>$ is a sequence of regular cardinals whose cofinalities are not eventually constant. 
 
 The first possibility is that some of them are not regular cardinals in $M^X_{\theta_X}$. There exists then some $M^* \propeextend M^X_{\theta_X}$ which projects below $\on \cap N^X_{\zeta_X}$. By minimizing we can have that all of $\<\gamma^X_{n_i + 1}: i < \omega\>$ are regular cardinals in $M^*$. As $M^*$ is sound and a tail of the $\gamma^X$'s lies above its projectum we can apply \thref{corelemma} to show that their cofinalities must eventually equal the cofinality of some fixed projectum of $M^*$. Contradiction!
 
 The second possibility is that the $\gamma^X$'s are regular in $M^X_{\theta_X}$. But then by assumption we have that $M^X_{\theta_X}$ is sound above some $\eta < \on \cap N^X_{\zeta_X}$ (the sup of generators of the iteration). This leads to the same contradiction as above!
\end{beweis}

So we can deduce that there must eventually be a final drop on the $K$-side. We will have that the $d^X_\cdot$ (, the least $m$ such that $\kappa^X_{\cdot} \geq \rho_{m+1}(M^X_{\cdot})$,) are constant on a tail-end of the iteration. We can assume w.l.o.g. that this eventual value is some fixed $m$ for all $X \in S$. Then $\cof(\rho_m(M^X_\cdot))$ also will have a stable end-value, call it $\rho$. 

There must then be some natural number $i^*$ such that $\cof(\rho_m(M^X_\alpha)) = \rho$ whenever $\nu^X_\alpha \geq \gamma^X_{n_{i^*}}$.

\begin{lemma}
 Let $i \geq i^*$ be such that $\mu_{n_i} \neq \rho$. Let $\alpha < \theta_X$ be minimal with $\nu^X_\alpha \geq \gamma^X_{n_i}$. Then $\kappa^X_\alpha < \gamma^X_{n_i}$.
\end{lemma}

\begin{beweis}
 Assume not. So $\kappa^X_\alpha \geq \gamma^X_{n_i}$. We then have that $\gamma^X_{n_i + 1}$ is a regular cardinal in $M^X_\alpha$ by the agreement between models in an iteration. Also $M^X_\alpha$ must be $(m+1)$-sound above $\gamma^X_{n_i}$ by the minimality of $\alpha$. Hence by \thref{corelemma} we have that $\cof(\gamma^X_{n_i + 1}) = \cof(\rho_m(M^X_\alpha)) = \rho$. This contradicts the choice of $i$!
\end{beweis}

One can even show that $((\kappa^X_\alpha)^+)^{M^X_\alpha} < \gamma^X_{n_i}$. One uses that if $\gamma^X_{n_i}$ is a succesor then its cofinality must equal $\mu_{n_i}$.

\begin{lemma}
 Let $i \geq i^*$ be such that $\mu_{n_i} \neq \rho$. Let $\alpha < \theta_X$ be minimal with $\nu^X_\alpha \geq \gamma^X_{n_i}$. Then $\nu^X_\alpha \geq \gamma^X_{n_i + (l - 1)}$.
\end{lemma}

\begin{beweis}
 Assume not. Then $((\nu^X_\alpha)^+)^{M^X_{\alpha + 1}}$ is a cardinal in $\left(\gamma^X_{n_i},\gamma^X_{n_i + (l-1)}\right]$ and thus has cofinality $\mu_{n_i} \neq \rho$. On the other hand $M^X_{\alpha + 1}$ is $(m+1)$-sound above $\nu^X_\alpha$ and $((\nu^X_\alpha)^+)^{M^X_{\alpha + 1}}$ is regular there. So then by \thref{corelemma} we have $\cof(((\nu^X_\alpha)^+)^{M^X_{\alpha + 1}}) = \cof(\rho_m(M^X_{\alpha + 1})) = \rho$. Contradiction!
\end{beweis}

So then a $\kappa^X_\alpha$ as above has Mitchell-order at least $((\kappa^X_\alpha)^{+(l+1)})^{N^X_{\zeta_X}}$. This then pulls back to $K_X$. Note that $\kappa^X_\alpha$ is in the range of $\tau^X_{0,\zeta_X}$ because it is the largest measurable cardinal below $\gamma^X_{n_i}$. 

Of course there are infinitely many $i$ such that $\mu_{n_i}$ is not equal to $\rho$. This finishes the proof of \thref{first}.

\thref{second} now follows rather quickly:

\begin{beweis}
	Fix natural numbers $k < l$. Assume that for all $f:\omega \rightarrow \{l,k\}$ the sequence $\<S^n_{f(n)} : l < n < \omega\>$ is mutually stationary. Assume for contradiction that $0^\P$ does not exist.
	
	The first thing we need to realize is that by the results of the last section there are unboundedly many $K$-measurable cardinals below $\aleph_\omega$.
	
	We can then find a sequence $\<I_n: n < \omega\>$ of ordinal intervals with cardinal endpoints that partition $\aleph_\omega$ such that there exists $\<\kappa_n: n < \omega\>$ $K$-measurable cardinals with $\{\kappa_n,(\kappa^+_n)^K\} \subset I_n$ for all $n < \omega$.
	
	Define a function $f:\omega \rightarrow \{l,k\}$ by 
	
	\begin{equation*}
		m \mapsto \begin{cases} l & m \in I_n \text{ where } n \text{ is odd} \\
                                k & m \in I_n \text{ where } n \text{ is even.} \end{cases}
    \end{equation*}
    
    By \thref{prop1} we can then find some $S \subset \Pot(H_{\aleph_\omega})$ stationary such that $S$ meets $\<S^n_{f(n)}: l < n < \omega\>$ on a tail, i.e. $\sup(X \cap \aleph_n) \in S^n_{f(n)}$ for all $X \in S$ and all but finitely many $n$, and for all $X \in S$ in the co-iteration between $K$ and $K_X$ the $K$-side truncates.
    
    Note that by choice of $f$ we have that $\cof(((\kappa^X_n)^+)^{K_X})$ is $\aleph_l$ iff $n$ is odd and $\aleph_k$ iff $n$ is even for all $X \in S$ and all but finitely many $n$ where $\kappa^X_n := \sigma^{-1}_X(\kappa_n)$.
    
    Let $(\kappa^X_n)^* := \tau^X_{0,\zeta_X}(\kappa^X_n)$.
    
    Again we can fix $m,n^*,\rho$ such that $\cof(\rho_m(M^X_\alpha)) = \rho$ whenever $\nu^X_\alpha \geq (\kappa^X_{n^*})^*$. 
    
    There is then some $q \geq n^*$ such that $\cof((((\kappa^X_q)^*)^+)^{N^X_{\zeta_X}}) \neq \rho$. Let then $\alpha < \theta_X$ be minimal with $\nu^X_\alpha \geq (\kappa^X_q)^*$. As before we have that $\kappa^X_\alpha < (\kappa^X_q)^*$. But then in $N^X_{\zeta_X}$ the Mitchell-order of $\kappa^X_\alpha$ is at least $(\kappa^X_q)^*$ which is measurable there. But then $0^\P$ exists. Contradiction!
\end{beweis}

\section{Very small J\'{o}nsson cardinals}

Now let us fix some J\'{o}nsson cardinal $\kappa$ with $\kappa < \aleph_k$ but assume that $0^\P$ does not exist. W.l.o.g. we do assume that $\kappa$ is a limit cardinal. Say $X \prec (H_\kappa; \in, K \cap \kappa)$ is a J\'{o}nsson type substructure, i.e. $\kappa \nsubseteq X$ but $\card(X \cap \kappa) = \kappa$. To simplify one argument in particular we do assume that the set of cardinals below $\kappa$ is contained in $X$.

We let $\mu^X_0$ be the least cardinal $\mu$ (smaller than $\kappa$) such that $X \cap \mu \neq \mu$(, alternatively $\mu^X_0 = \sigma_X(\crit(\sigma_X))$). Let $\mu^X_{1 + \alpha}$ be the least cardinal $\mu$ such that $X \cap \mu$ has size $(\mu^X_0)^{+\alpha}$ for $\alpha < \alpha^*$ where $\alpha^*$ is such that $(\mu^X_0)^{+\alpha^*} = \kappa$.

\begin{proposition}
	$\mu^X_{(\alpha + 1)}$ is a successor cardinal for all but boundedly many $\alpha < \alpha^*$.
\end{proposition}

\begin{beweis}
	 Assume not. Say $\lambda := \mu^X_{(\alpha + 1)}$ is a limit cardinal, then we must have $\cof(X \cap \lambda) = (\mu^X_0)^{+\alpha}$. Therefore $\cof(\lambda) = (\mu^X_0)^{+\alpha}$ and hence $\lambda \geq \aleph_{(\mu^X_0)^{+\alpha}}$. As the sequence $\<(\mu^X_0)^{+\alpha}:\alpha < \alpha^*\>$ converges to $\kappa$ we have $\kappa = \aleph_\kappa$. Contradiction!
\end{beweis}

Our proof requires only that there exists a sequence $\<\alpha_n : n < \omega\>$ such that $\mu^X_{(\alpha_n + 1)}$ is a successor cardinal for all $n < \omega$. We will assume that $\alpha_n = n$, the general case is only notationally more complex. 

We will obviously have that $\cof(X \cap \mu^X_{n+1}) = (\mu^X_0)^{+n}$ but it is also true that $\cof(X \cap (\mu^X_{n+1})^{-}) = (\mu^X_0)^{+(n-1)}$ where $(\mu^X_{n+1})^{-}$ is the cardinal predecessor of $\mu^X_{n+1}$ assuming that predecessor is regular.

The second is immediate if $(\mu^X_{n+2})^- = \mu^X_{n+1}$ but is true for less immediate reasons otherwise by work of Shelah \cite[p 451]{Jech} adapted to our current context by Foreman and Magidor \cite[2.15]{fandmdefcounterexampletoch}.

As usual our proof is based on analysing the co-iteration between $K$ and $K_X$. We will first assume that $K$ early on in the iteration truncates to a mouse of size $\kleiner\mu^X_0$. We do not require that the $K_X$ side remains trivial.

We will assume that $\zeta_X = \kappa$, otherwise we can always pad the iteration. As $K$ "wins" the iteration but does by assumption truncate to a small mouse we must have clubs $C^X_n \subset (\mu^X_0)^{+n}$ such that $\kappa^X_\alpha = \alpha$, $\pi^X_{\alpha,(\mu^X_0)^{+n}}(\alpha) = (\mu^X_0)^{+n}$, $d^X_\alpha$ is constant and no truncation takes place in the interval $\left[\alpha,(\mu^X_0)^{+n}\right]$ for all $\alpha \in C^X_n$. (We must have $(\mu^X_0)^{+n} \in M^X_{(\mu^X_0)^{+n}}$. If $\gamma$ then is some pre-image of $(\mu^X_0)^{+n}$ in $M^X_\beta$, then an appropriate subset of $\{\pi^X_{\beta,\alpha}(\gamma) \vert \beta \leq \gamma < (\mu^X_0)^{+n}\}$ is as wanted.)

By the proof of the comparison lemma we cannot have simultaneously some $\beta,\gamma < (\mu^X_0)^{+n}$ such that $\tau^X_{\beta,(\mu^X_0)^{+n}}(\gamma) = (\mu^X_0)^{+n}$. Therefore we have clubs $D^X_n \subset (\mu^X_0)^{+n}$ such that $\tau^X_{0,\alpha}(\alpha) = \alpha$ for all $\alpha \in D^X_n$. (We will also assume that $D^X_n \subseteq \left(\sigma^{-1}_X((\mu^X_{n+1})^-),(\mu^X_0)^{+n}\right)$.)

We can then immediately conclude that $(\mu^X_{n+1})^{-}$ is not singular as otherwise $\mu^X_{n+1}$ would be a successor cardinal in $K_X$ but every element in $C^X_{n+1} \cap D^X_{n+1}$ is a $K_X$-cardinal. Therefore $\cof(\sigma^{-1}_X((\mu^X_{n+1})^{-})) = (\mu^X_0)^{+(n-1)}$ for all $n < \omega$.

\begin{lemma}
	Let $n < \omega$:
	\begin{itemize}
		\item[$(a)$] $\forall \alpha \in D^X_n : \cof((\alpha^+)^{N^X_\alpha}) = (\mu^X_0)^{+(n-1)}$;
		\item[$(b)$] $\forall \alpha < (\mu^X_0)^{+n} : \cof((\alpha^+)^{M^X_\alpha}) = \cof(\rho_{d^X_\alpha}(M^X_\alpha))$.
	\end{itemize}
\end{lemma}

\begin{beweis}
	
	\begin{itemize}
		\item[$(a)$] It follows from weak covering that $\cof((\alpha^+)^{K_X}) = (\mu^X_0)^{+(n-1)}$. $\tau_{0,\alpha}$ is continuous at $(\alpha^+)^{K_X}$ and $\alpha$ itself is fixed by the embedding, so $\cof((\alpha^+)^{N^X_\alpha}) = (\mu^X_0)^{+(n-1)}$.
		\item[$(b)$] This follows from applying \thref{corelemma} in $M^X_\alpha$ keeping in mind that the generators of the iteration up to stage $\alpha$ are bounded by $\kappa^X_\alpha = \alpha$.
	\end{itemize}
	
\end{beweis}

An immediate corollary to the lemma is that for all $n < \omega$, all $\alpha \in D^X_n \cap C^X_n$ we have $\cof(\rho_{d^X_\alpha}(M^X_\alpha)) = (\mu^X_0)^{+(n-1)}$. ($\alpha \in C^X_n$ implies that we do not truncate when forming $M^X_{\alpha + 1}$, so $(\alpha^+)^{M^X_\alpha} = (\alpha^+)^{N^X_\alpha}$.)

We conclude that the iteration map between $\alpha \in C^X_n \cap D^X_n$ and $\beta \in C^X_{n+1} \cap D^X_{n+1}$ did not move $\rho_{d^X_\alpha}(M^X_\alpha)$ cofinally to $\rho_{d^X_\beta}(M^X_\beta)$. Therefore we dropped to a smaller model (either in model or degree) somewhere in the interval $\left(\alpha,\beta\right]$. But there can only be a finite amount of such drops in the course of an iteration. Contradiction!

We must therefore have that $K$ does not truncate to a model of size $\kleiner\mu^X_0$ in the course of the co-iteration. Following the proof of \cite[Lemma 39]{Cox} this can only be because of one of two reasons:

\begin{itemize}
	\item $\alpha^X_0 := (\crit(\sigma_X)^+)^{K_X} = (\crit(\sigma_X)^+)^K$ but $\ult(K;\sigma_X \restr (K_X \vert\vert \alpha^X_0))$ is not iterable;
	\item $\alpha^X_0 < (\crit(\sigma_X)^+)^K$, $\alpha^X_1 := ((\nu^X_0)^+)^{K_X} = ((\nu^X_0)^+)^{M^X_1}$ but $\ult(M^X_1;\sigma_X \restr (K_X \vert\vert \alpha^X_1))$ is not iterable.
\end{itemize}

We can effectively ignore the first case: w.l.o.g. we can assume that $\mu^X_0 \geq \aleph_2$; then $\cof(\alpha^X_0) \geq \aleph_1$; it follows from \cite[Corollary 19]{Cox} that $\ult(K;\sigma_X \restr (K_X \vert\vert \alpha^X_0))$ must then be iterable. 

For the second case we shall apply an idea from \cite{mitchelljonsson}: fix some $\<f_n,a_n: n < \omega\>$ witnessing the fact that $\ult(M^X_1;\sigma_X \restr (K_X \vert\vert \alpha^X_1))$ is not iterable. Let $\gamma$ be a regular cardinal with $\<f_n:n < \omega\> \subset M^X_1\vert\vert\gamma$.
 
We let $M^*$ be the transitive collapse of $\hull^{M^X_1 \vert\vert \gamma}_\omega(\<f_n:n < \omega\> \cup \alpha^X_1\>$. We have that $\card(M^*) < \mu^X_0$ and $\ult(M^*; \sigma_X \restr \alpha^X_1)$ is not iterable as witnessed by $\<\bar{f}_n,a_n: n < \omega\>$ where $\bar{f}_n$ is the image of $f_n$ under the collapse.
 
\begin{lemma}
 $M^*$ wins the co-iteration with $K_X$.
\end{lemma}

\begin{beweis}
Assume not. Let $\<N_\alpha, \tau_{\alpha,\beta}: \alpha \leq \beta \leq \zeta^*\>$ be the iteration on $K_X$ and $\<M^*_\alpha, \pi^*_{\alpha,\beta}: \alpha \leq \beta \leq \theta^*\>$ be the iteration on $M^*$. By assumption we have that $M^*_{\theta^*} \eextend N_{\zeta^*}$.

Let $\<N^*_\alpha,\tau^*_{\alpha,\beta}:\alpha \leq \beta \leq \zeta^*\>$ be the iteration on $K$ that results from copying the iteration of $K_X$ by $\sigma_X$. Let $\<\sigma^*_\alpha: \alpha \leq \zeta^*\>$ be the copy maps.

Clearly, $\ult(N_{\zeta^*},\sigma_{\zeta^*} \restr (K_X \vert\vert \alpha^X_1))$ is iterable and thus so is $\ult(M^*_{\theta^*},\sigma_{\zeta^*} \restr (K_X \vert\vert \alpha^X_1))$. 

Note that $\nu^X_0$ is a strong cutpoint of $K_X$ and $M^X_1$ so the iteration is above $\alpha^X_1$. Therefore we can embed $\ult(M^*;\sigma_X \restr (K_X \vert\vert \alpha^X_1))$ into $\ult(M^*_{\theta^*}; \sigma^*_{\zeta^*} \restr (K_X \vert\vert \alpha^X_1))$ ($(f,a) \mapsto (\pi^*_{0,\theta^*}(f),a)$). But by choice of $M^*$ the former is \textit{not} iterable. Contradiction!
\end{beweis}

From this point on we need not concern ourselves with $K$ any longer, so we will re-index and use $\<N_\alpha,\tau_{\alpha,\beta}: \alpha \leq \beta \leq \zeta_X\>$ and $\<M^X_i,\kappa^X_i,\nu^X_i,\pi^X_{i,j},d^X_i: i \leq j \leq \theta_X\>$ for the iterations on $K_X$ and $M^*$ respectively that arise out of the co-iteration of the two models. (We will remember the first two steps of the previous iteration, so $M^* = M^X_2$.)

As the $M^*$-side is small we can find once again clubs $C^X_n \subset (\mu^X_0)^{+n}$ consisting of the preimages of $(\mu^X_0)^{+n}$ under $\pi^X_\cdot$ and $D^X_n \subset (\mu^X_0)^{+n}$ consisting of fixpoints of $\tau^X_\cdot$. 

As before we quickly realize that for no $n < \omega$ is $\mu^X_{n+1}$ the successor of a singular cardinal. We then have that $\cof(\sigma^{-1}_X((\mu^X_{n+1})^{-})) = (\mu^X_0)^{+(n-1)}$ for all $n < \omega$. As before we conclude that $\cof((\alpha^+)^{N^X_\alpha}) = (\mu^X_0)^{+(n-1)}$ for $\alpha \in D^X_n$.

Mitchell characterizes the move from $M^X_1$ to $M^*$ as a "drop" but unfortunately $M^*$ does lack some characteristics typical for drops, e.g. $\rho_\omega(M^*) = \on \cap M^*$. This means that for once our usually so reliable \thref{corelemma} deserts us. Fortunately, $M^*$ does share the characteristic of being a minimal counter-example of sorts that we can exploit to prove a substitute lemma.

\begin{lemma}
	Let $\alpha < \theta_X$ and assume that no drop occurs in the interval $\left[2,\alpha\right)$, then $\cof((\kappa^+_\alpha)^{M^X_\alpha}) \leq \cof(\on \cap M^X_\alpha)$.
\end{lemma}

\begin{beweis}
	Fix some such $\alpha$. For $\beta < \on \cap M^X_\alpha$ we define $\gamma^n_\beta := \sup(\hull^{M^X_\alpha \vert\vert\beta}_\omega(\kappa_\alpha \cup \{\pi^X_{2,\alpha}(\bar{f}_n):m < n\}) \cap (\kappa^+_\alpha)^{M^X_\alpha})$. If there were to be some $\beta$ such that $\gamma_\beta := \sup\limits_{n < \omega} \gamma^n_\beta = (\kappa^+_\alpha)^{M^X_\alpha}$ then the latter would have cofinality $\omega \leq \cof(\on \cap M^X_\alpha)$. So assume $\gamma_\beta < (\kappa^+_\alpha)^{M^X_\alpha}$ for all $\beta < \on \cap M^X_\alpha$.
	
	It will be sufficient to prove that $\<\gamma_\beta: \beta < \on \cap M^X_\alpha\>$ is cofinal in $(\kappa^+_\alpha)^{M^X_\alpha}$. So fix $\delta < (\kappa^+_\alpha)^{M^X_\alpha}$. An easy induction shows that $M^X_\alpha = \hull^{M^X_\alpha}_\omega(\kappa^X_\alpha \cup \<\pi_{2,\alpha}(\bar{f}_n):n < \omega\>)$. Fix then some first order formula $\phi$ together with $\xi_0,\ldots,\xi_{m-1}$ such that $\delta$ is the least ordinal with $\phi^{M^X_\alpha}(\delta,\xi_0,\ldots,\xi_{m-1},\pi_{2,\alpha}(\bar{f}_0),\ldots,\pi_{2,\alpha}(\bar{f}_{m-1}))$.
	
	We claim that 
	\begin{align*}
	\forall \delta \forall \xi_0 \ldots \forall \xi_{m-1} \exists \beta: & \phi^{M^X_\alpha \vert\vert \beta}(\delta,\xi_0,\ldots,\xi_{m-1},\pi_{2,\alpha}(\bar{f}_0),\ldots,\pi_{2,\alpha}(\bar{f}_{m-1})) \\
	                                                                \gdw & \phi(\delta,\xi_0,\ldots,\xi_{m-1},\pi_{2,\alpha}(\bar{f}_0),\ldots,\pi_{2,\alpha}(\bar{f}_{m-1})) \end{align*}
	holds in $M^X_\alpha$. This is a first order statement that holds in $M^X_1 \vert\vert \gamma$ (relative to $f_0,\ldots,f_{n-1}$) as $\gamma$ was a regular cardinal. Thus it also holds in $M^*$ relative to $\bar{f}_0,\ldots,\bar{f}_{n-1}$. By assumption $\pi_{2,\alpha}$ is fully elementary so it must also hold in $M^X_\alpha$. Take then some such $\beta$ and we will have $\delta < \gamma_\beta$ as desired.
\end{beweis}

As $M^*$ has size $\kleiner\mu^X_0$ we then consequently also must have $\cof(\on \cap M^*) < \mu^X_0$. The latter fact will be preserved by the embeddings $\pi_{2,\alpha}$ if they exist. We can thus conclude that the iteration truncates somewhere before $\min(C^X_1 \cap D^X_1)$. But past this truncation our (as Mitchell calls it) "quasi-iteration" does not meaningfully differ from usual iterations. So we can argue just as before that there must be infinitely more truncations. But that is a contradiction! 

\section{Open Questions}

\begin{frage}
	Assume that the hypothesis of \thref{first} hold for some $2 \leq l < \omega$: is $\mitord^K(\aleph_\omega) > (\aleph^{+(l+1)}_\omega)^K$?
\end{frage}

\begin{frage}
	Let $1 \leq k < l < \omega$. Is it consistent relative to large cardinals that for all $f:\omega \rightarrow \{k,l\}$ the sequence $\<S^n_{f(n)}: k < n < \omega\>$ is mutually stationary?
\end{frage}

\begin{frage}
	Assume that $\kappa$ is a J\'{o}nsson cardinal with $\kappa < \aleph_\kappa$. Does there exist an inner model with a Woodin cardinal?
\end{frage}

\section{Acknowledgements and historical note}

The first version of \thref{first} emerged while the author was working at UC Berkeley under the DFG grant AD 469/1-1 (while work on what would become \cite{altmutstat} was still on-going). In this earlier version we could show that there must be some $\lambda$ in $K$ that was a limit of cardinals $\kappa$ with $\mitord(\kappa) \geq \kappa^{+(l-1)}$, but $\lambda < \aleph_\omega$ seemed possible. Therefore this version did not lead to a version of \thref{second}.

The current version of \thref{first} (and with it \thref{second}) came to the author while he was a participant in the program "Large Cardinals and Strong logics" taking place in and around the CRM in Bellaterra, Spain during the winter of 2016. 

This paper was finally finished while the author was a postdoc at Bar-Ilan University in Ramat Gan, Israel with a grant from the European Research Council (grant agreement ERC-2018-StG 802756).

\bibliographystyle{alpha}
\bibliography{bibliographie}

\begin{thebibliography}{CFM06}

\bibitem[ACW]{altmutstat}
Dominik Adolf, Sean Cox, and Philip Welch.
\newblock Lower consistency bounds for mutual stationarity with divergent
  uncountable cofinalities.
\newblock to appear, online preprint: https://arxiv.org/abs/1607.04790.

\bibitem[Bau91]{thatonelemma}
J.E. Baumgartner.
\newblock On the size of closed unbounded sets.
\newblock {\em Annals of Pure and Applied Logic}, 54:195 -- 227, 1991.

\bibitem[BN]{singularstatI}
Omer Ben-Neria.
\newblock On singular stationarity i.
\newblock submitted.

\bibitem[CFM06]{canstrtwo}
James Cummings, Matthew Foreman, and Menachem Magidor.
\newblock Canonical structure in the universe of set theory: part two.
\newblock {\em Annals of Pure and Applied Logic}, 142:55 -- 75, 2006.

\bibitem[Cox09]{Cox}
Sean Cox.
\newblock Covering theorems for the core model, and an application to
  stationary set reflection.
\newblock {\em Annals of Pure and Applied Logic}, 161:66 -- 93, 2009.

\bibitem[FM95]{fandmdefcounterexampletoch}
Matthew Foreman and Menachem Magidor.
\newblock Large cardinals and definable counterexamples to the continuum
  hypothesis.
\newblock {\em Annals of Pure and Applied Logic}, 76(1):47 -- 97, 1995.

\bibitem[FM01]{mutstat}
Matthew Foreman and Menachem Magidor.
\newblock Mutually stationary sequences of sets and the non-saturation of the
  non-stationary ideal on {$\Pot_\kappa(\lambda)$}.
\newblock {\em Acta Mathematica}, 186(2):271 -- 300, 2001.

\bibitem[Jec06]{Jech}
Thomas Jech.
\newblock {\em Set Theory - The Third Millenium Edition,revised and expanded}.
\newblock Springer Monographs in Mathematics. Springer-Verlag, Berlin, 3rd
  edition, 2006.

\bibitem[Jen]{chforcax}
Ronald Jensen.
\newblock Forcing axioms compatible with {$\CH$}.
\newblock online at http://www.mathematik.hu-berlin.de/~raesch/org/jensen.html.

\bibitem[KW06]{strengthmutstat}
Peter Koepke and Philip Welch.
\newblock On the strength of mutual stationarity.
\newblock In {\em Set Theory: Centre de Recerca Matematica Barcelona,
  2003--2004}, Trends in Mathematics. Birkh\"{a}user, 2006.

\bibitem[LS97]{liushelah}
Kecheng Liu and Saharon Shelah.
\newblock Cofinalities of elementary substructures of structures on
  $\aleph_\omega$.
\newblock {\em Israel Journal of Mathematics}, 99(1):189 -- 205, 1997.

\bibitem[Mit99]{mitchelljonsson}
William Mitchell.
\newblock J\'{o}nsson cardinals, {Erd\H{o}s} cardinals, and the {Core Model}.
\newblock {\em The Journal of Symbolic Logic}, 64(3):1065 -- 1086, 1999.

\bibitem[Sha00]{sharonmsth}
Assaf Sharon.
\newblock Generators of {$\PCF$}.
\newblock Master's thesis, Tel Aviv University, 2000.

\bibitem[Zem01]{Zeman}
Martin Zeman.
\newblock {\em Inner Models and large Cardinals}, volume~5 of {\em De Gruyter
  Series in Logic and its applications}.
\newblock De Gruyter, Berlin;New York, 2001.

\end{thebibliography}

\end{document}